\documentclass[a4paper,12pt]{article}
\usepackage[top=2.5cm,bottom=2.5cm,left=2.5cm,right=2.5cm]{geometry}
\usepackage{cite, amsmath, amssymb}
\usepackage[margin=1cm,font=small,format=hang,labelsep=period,labelfont=bf]{caption}
\usepackage{tikz}
\usepackage{authblk}
\usepackage{footnote}
\usepackage{amsmath}
\usepackage{amsfonts}
\usepackage{amssymb}
\usepackage{graphicx}
\usepackage{enumitem}

\newtheorem{thr}{Theorem}[section]

\newtheorem{prop}[thr]{Proposition}
\newtheorem{defi}[thr]{Definition}

\def\U{\mathbb{U}}
\def\T{\mathbb{T}}
\def\P{\mathcal{P}}

\title{Maximum Wiener Indices of Unicyclic Graphs
	of Given Matching Number}
\author{Stijn Cambie \footnote{This work has been supported by a Vidi Grant of the Netherlands Organization for Scientific Research\\ (NWO), grant number $639.032.614$ }}%
\affil{Department of Mathematics, Radboud University Nijmegen, Postbus 9010, 6500 GL Nijmegen, The Netherlands; S.Cambie@math.ru.nl.}
\date{}
\begin{document}

	%\nonumber
\tikzstyle{every node}=[circle, draw, fill=black!50,
	inner sep=0pt, minimum width=4pt]

%\bigskip[6 cm]
%\hspace{5cm}	
%\mbox{}
%\vspace{160pt}

	\maketitle
	\begin{abstract}
		In this article, we determine the maximum Wiener indices
		of unicyclic graphs with given number of vertices and matching number.
		We also characterize the extremal graphs. % attaining it.
		%This solves the remaining open problem concerning Wiener 
		This solves an open problem of Du and Zhou~\cite{DZ}.
	\end{abstract}

    	\section{Introduction}
Let $G$ be a simple connected graph. We denote its vertex set by $V(G)$ and its edge set by $E(G).$
For two graphs $G$ and $H$, $G+H$ denotes the graph obtained by adding all possible edges between the vertices of $G$ and $H$.
For an integer $p$, $pG$ is the disjoint union of $p$ copies of $G.$
The matching number of a graph $G$ is the size of a maximum independent edge subset of $G$, we will denote it by $m(G)$ or $m$.
We will denote by $\T(n,m)$ and $\U(n,m)$ the set of all trees and unicyclic graphs respectively with $n$ vertices and matching number $m$. 
Let $d(u, v)$ denote the distance between vertices $u$ and $v$ in a graph $G$.
The Wiener index of a graph $G$ equals the sum of distances between all unordered pairs of vertices, i.e. $$W(G)=\sum_{\{u,v\} \subset V(G)} d(u,v).$$
The Wiener index, introduced in $1947$ by Harry Wiener, is the oldest and one of the most important topological indices in chemical graph theory. This index can predict some chemical properties of molecules, e.g. the boiling point of alkanes, the density, the critical point and surface tension. It was used by chemists decades before it attracted the attention of mathematicians.
An overview of the mathematical results, conjectures and problems with respect to the Wiener index can be found in \cite{KST}. 
In chemical graph theory, an important goal is to bound some important graph invariant like the Wiener index using some structural parameters. 
In this paper we consider the case of the matching number and we solve the remaining open case, Problem 11.4 of \cite{KST}. 

The minimum and maximum value for the Wiener index of a connected graph of order $n$ with matching number $m$ was determined by Dankelmann~\cite{D}.
% The result is the following.

\begin{thr}[Dankelmann~\cite{D}]\label{ThrDankelmann}
	Let $G$ be a connected graph of order $n$ and matching number $m$.
	\begin{itemize}
		\item 	If $m= \lfloor \frac n2 \rfloor$, then $W(G) \ge \binom{n}{2}$. Equality holds iff $G$ is the complete graph $K_n$.
		\item  	If $1 \le m < \lfloor \frac n2 \rfloor$, then $W(G) \ge 2\binom{n}{2}-mn+\binom{m}{2}$. Equality holds iff $G$ is isomorphic to the graph obtained by $K_m+(n-m)K_1$.
		\item   If $n$ is even, then $W(G) \le \binom{2m}{3}+\binom{2m}2 (n-2m+1)+2m\frac{n-2m+2}{2}\frac{n-2m}{2}+\frac12(n-2m)^2$.  Equality holds iff the graph is $A_{n,m}$, a path with $2m-1$ vertices, with one leaf of the path connected to $\frac{n-2m+2}{2}$ vertices and the other leaf with $\frac{n-2m}{2}$ vertices.
		\item   If $n$ is odd, then $W(G) \le \binom{2m}{3}+\binom{2m}2 (n-2m+1) + 2m(\frac{n-2m+1}{2})^2+4\binom{\frac{n-2m+1}{2}}{2}$.
		Equality holds iff the graph is $A_{n,m}$, a path with $2m-1$ vertices, with both leafs of the path connected to $\frac{n-2m+1}{2}$ different vertices, equivalently a path with $2m-3$ vertices with its ends concatenated to the centers of two stars of order $\frac{n-2m+3}{2}$.
	\end{itemize}
	
\end{thr}

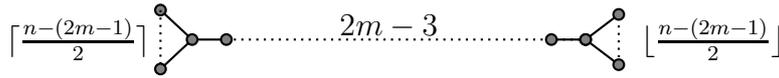
\begin{figure}[h]
	\centering
	\begin{tikzpicture}[thick,scale=0.9]%
	\draw \foreach \a in {10,-10} \foreach \x in {180} {
		(\x:2)--(\x+\a:2.5) node {}
	};
	\draw[dotted] \foreach \x in {0} {
		(-1.5,0)--(3.75,0)  
	};
	\draw 	(0:3.25)--  (0:3.75)  node {};
	\draw 	(-2,0)--  (-1.5,0)  node {};
	\draw[dotted] (170:2.5)--(190:2.5);

	\draw[dotted]  \foreach \x in {0} {
		(\x-5:4.25)--(\x+5:4.25) 
	};
	
	\draw \foreach \a in {5,-5} \foreach \x in {0} {
		(\x:3.75)--(\x+\a:4.25) node {}
	};
	
	%		\coordinate [label=left:$A_0$] (A) at (2.5,0.35); 
	%		\coordinate [label=left:$B_0$] (A) at (0.85,2.1); 
	%		\coordinate [label=right:$C_0$] (A) at (-1.96,0); 
	%		\coordinate [label=left:$D_0$] (A) at (0.8,-2.1); 
	%		\coordinate [label=left:$A_1$] (A) at (3.15,0.35); 
	%		
	\coordinate [label=center:$2m-3$] (A) at (0.875,0.2);

	\coordinate [label=right:$\lfloor \frac{n-(2m-1)}{2} \rfloor $] (A) at (0:4.5); 
	\coordinate [label=left:$\lceil \frac{n-(2m-1)}{2} \rceil $] (A) at (180:2.55); 
	%		\draw \foreach \x in {0,90,180,270} {
	%			%(\x:3) node {} -- (\x+90,3))
	%			(\x:2) node {}  -- (\x+90:2)
	%			
	%		};
	
	\draw	(3.25,0)  node {};
	\draw	(-2,0)  node {};	
	\draw	(3.75,0)  node {};

	\end{tikzpicture}
	\caption{Extremal graph $A_{n,m}$}
	\label{fig:graphDankelmann}
\end{figure}

The graphs that attain the maximal Wiener indices are trees, while the graphs that attain the minimal Wiener indices contain many cycles.
Hence two natural questions arise. What are the minimal Wiener indices for trees? 
What are the maximal Wiener indices for the graphs that are not trees? In the latter case, the graphs will be unicyclic. 
For graphs that contain multiple cycles, we can take a maximal matching, and then deleting an edge not contained in the matching such that the graph is still connected will increase the Wiener index.

Du and Zhou~\cite{DZ} solved the first question and also determined the minimal matching number for unicyclic graphs.
Their results are stated in the following two theorems.

\begin{thr}[Du and Zhou~\cite{DZ}]
	Let $G\in \T(n,m)$, where $2 \le m \le \lfloor \frac n2 \rfloor$.
	Then $W(G) \ge n^2 + (m-3)n-3m+4.$
	Equality holds iff $G$ is a star of order $n-m+1$ with an additional vertex connected to each of $m-1$ leaves.
\end{thr}

\begin{thr}[Du and Zhou~\cite{DZ}]
	Let $G\in \U(n,m)$, where $2 \le m \le \lfloor \frac n2 \rfloor$.
	If $(n,m)=(6,3)$, then $W(G) \ge 26,$ with equality iff $G$ is $C_5$ with a pendent vertex attached to it.
	In the other cases $W(G) \ge n^2+(m-4)n-3m+6.$
	Equality holds for a star of order $n-m$ with a triangle attached to its center and an additional vertex connected to each of $m-2$ leaves, as well as for $C_4, C_5$, the graph $C_5$ with two pendent vertices attached to it, and $C_5$ with three pendent vertices attached to three consecutive vertices of the cycle.
\end{thr}

In this paper, we will prove the remaining question, which was left as an open problem in \cite{DZ}.
The result is summarized in the following theorem.

\begin{thr}
	Let $G \in \U(n,m)$, where $2 \le m \le \lfloor \frac n2 \rfloor$.
	
	\begin{itemize}
		\item If  $n \le 2m+2$, then $W(G) \le 2-\frac83 m^3+2m^2+\frac53 m+2m^2n-3mn-2n+n^2$.
		\item If $n \ge 2m+3$ and $n$ is odd, then $W(G) \le \frac92-n-\frac23 m^3+\frac12 n^2-2nm+2m^2+\frac12mn^2-\frac{11}6m.$
		\item If $n \ge 2m+4$ and $n$ is even, then $W(G) \le 6-n-\frac23 m^3+\frac12 n^2-2nm+2m^2+\frac12 mn^2-\frac73m.$
	\end{itemize}
	Moreover, these bounds are sharp.	
\end{thr}

\begin{figure}[h]
	\centering
	
	\begin{tikzpicture}[thick,scale=0.9]%
	\draw[dotted] 	\foreach \x in {0} {
		(\x:2.75)--(\x:3.75) 
	};
	\draw[dotted]  \foreach \x in {0} {
		(\x-6:4.25)--(\x+6:4.25) 
	};
	\draw \foreach \x in {0} {
		%(\x:3) node {} -- (\x+120,3))
		
		(\x:2)--(\x:2.75)  node {}
	};
	\draw (120:2)   -- (0:2) node {};
	\draw (240:2)   -- (120:2) node {};
	\draw (0:2)   -- (240:2) node {};
	\draw (0:2) node {};

	\draw \foreach \a in {6,-6} \foreach \x in {0} {
		(\x:3.75)--(\x+\a:4.25) node {}
	};

	\draw \foreach \x in {0} {
		%\draw[dotted] (\x:2.75)--(\x:3.5)  
		(\x:3.75)  node {}
		
	};

	\coordinate [label=left:$A_0$] (A) at (10.1:2.3); 
	\coordinate [label=right:$B_0$] (A) at (117.5:2.15); 
	\coordinate [label=right:$C_0$] (A) at (242.5:2.15); 
	\coordinate [label=left:$A_1$] (A) at (6.75:3); 
	%	\coordinate [label=right:$B_1$] (A) at (117.5:2.8); 
	%	\coordinate [label=right:$C_1$] (A) at (242.5:2.8); 
	\coordinate [label=left:$A_j$] (A) at (5.5:3.95); 
	%	\coordinate [label=right:$B_k$] (A) at (117.5:3.7); 
	%	\coordinate [label=right:$C_l$] (A) at (242.5:3.7); 

	\coordinate [label=center:$a$] (A) at (0:4.5);

	\end{tikzpicture} \quad
	\centering
	\begin{tikzpicture}[thick,scale=0.9]%

	\draw[dotted] \foreach \x in {0} {
		(\x:2.75)--(\x:3.75)  
	};
	\draw 	(0:2)--  (0:2.75)  node {};
	\draw[dotted] (170:2.5)--(190:2.5);
	
	\draw \foreach \a in {10,-10} \foreach \x in {180} {
		(\x:2)--(\x+\a:2.5) node {}
	};

	\draw[dotted]  \foreach \x in {0} {
		(\x-5:4.25)--(\x+5:4.25) 
	};
	
	\draw \foreach \a in {5,-5} \foreach \x in {0} {
		(\x:3.75)--(\x+\a:4.25) node {}
	};
	
	\coordinate [label=left:$A_0$] (A) at (2.5,0.35); 
	\coordinate [label=left:$B_0$] (A) at (0.85,2.1); 
	\coordinate [label=right:$C_0$] (A) at (-1.96,0); 
	\coordinate [label=left:$D_0$] (A) at (0.8,-2.1); 
	\coordinate [label=left:$A_1$] (A) at (3.15,0.35); 
	
	\coordinate [label=left:$A_j$] (A) at (3.85,0.3);

	\coordinate [label=center:$a$] (A) at (0:4.5); 
	\coordinate [label=left:$c$] (A) at (180:2.55); 
	\draw \foreach \x in {0,90,180,270} {
		%(\x:3) node {} -- (\x+90,3))
		(\x:2) node {}  -- (\x+90:2)
		
	};
	\draw \foreach \x in {0} {
		%\draw[dotted] (\x:2.75)--(\x:3.5)  
		(\x:3.75)  node {}
		
	};
	\end{tikzpicture}

	\caption{The graphs $G^3_{a,j}$ and $G^4_{a,c,j}$}
	\label{fig:G3aj_G4acj}
\end{figure}
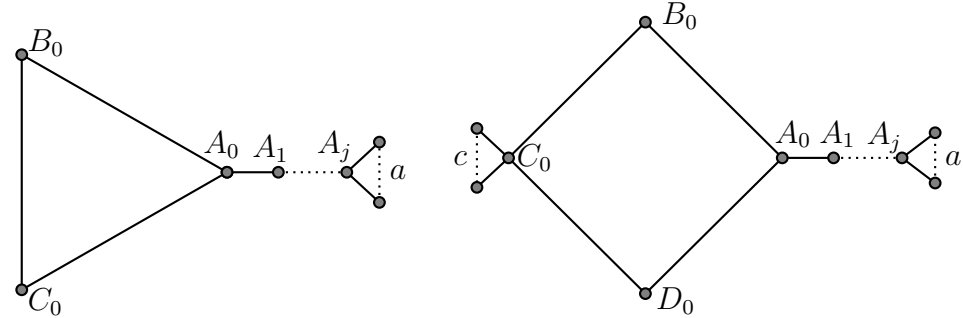

The extremal graphs are members of two families of graphs, which are shown in Figure~\ref{fig:G3aj_G4acj}. 
We describe them explicitly  in Theorem~\ref{mainthr}.

So in this paper, we determine the maximum Wiener indices of graphs in $\U(n,m)$, where $m\ge 2$ and $n \ge 2m$ and characterize the extremal graphs.
Note that the only nonempty remaining case is $(n,m)=(3,1)$, but then $C_3$ is the only unicyclic graph with these parameters. We will use the notation $W(\U(n,m)) = \max \{ W(G) \mid G \in \U(n,m)\}$.

Our proof for the characterization of $W(\U(n,m))$ will actually first determine the characteristics of the extremal graphs and only after that, we calculate the sharp upper bound.
It contains five main parts, which can be summarized as follows.
\begin{enumerate}
	\item We prove that for fixed $n$, $W(\U(n,m))$ increases when $m$ increases, as long as $m \le \lfloor \frac n2 \rfloor$.
	Intuitively, this is because a larger matching number implies that it is possible there are longer paths in the graph. 
	A path is the graph with the largest Wiener index for a fixed order.
	The proof for this property is given in Section~\ref{2}.
	\item Second, we prove that the cycle in the extremal graph in $\U(n,m)$ for any $n$ and $m$ must be small. More precisely this cycle must be a triangle $C_3$ or quadrangle $C_4$.
	The possible constructions for this are given in Section~\ref{3}. 
	We need two possible constructions to ensure that we do not increase the matching number.
	\item A unicyclic graph is a cycle and some attached trees. In Section~\ref{4} we prove that extremal graphs have attached trees which are of a certain form: they are the concatenation of a path and a star.
	Intuitively, this type of tree will have their vertices as far as possible from the remaining cycle and other attached trees. In this way the distances to the vertices outside the tree are as large as possible.
	\item Using the two previous steps, we know that any extremal graph can be represented with only a few parameters. In Section~\ref{5} and Section~\ref{6} we determine some conditions on these parameters, which are the length and the number of leaves of the subtrees.
	\item Finally, in Section~\ref{7} we calculate the Wiener index for the possible extremal graphs in $\U(n,m)$ and compare them to find the maximal value $W(\U(n,m)).$
\end{enumerate}

In the first three main parts, we are proving some heuristic arguments. To prove some of these arguments, we use tree rearrangements. In particular, we are using some kind of subtree pruning and regrafting (SPR), which we define here.
%https://en.wikipedia.org/wiki/Tree_rearrangement

\begin{defi}[SPR]
	Let $G$ be a graph.
	Given a rooted subtree $S$ of $G$, such that the root $d=S \cap H.$
	Pruning $S$ from $G$ is removing the whole structure $S$ excluding the root $d$.
	Regrafting $S$ at a vertex $v$, means that we are taking a copy $S'$ of $S$ which we insert at $v$, letting its root $d'$ coincide with $v$. No additional edges are drawn in this process.
\end{defi}

\begin{figure}[h]
	\centering
	
	\begin{tikzpicture}[thick,scale=0.9]%
	\draw[dotted] (1.2,1.2)--(1,0)--(0.8,1.2);
	\draw[dashed] (-1,0) -- (0,0);
	\draw[dashed] (0,0) -- (1,0);
	\draw[dashed] (1,0)  -- (2,0);
	\draw node{} (0,0) node{} ;
	\draw node{} (1,0) node{} ;
	\coordinate [label=center:$S$] (A) at (1,1.4);
	\coordinate [label=center:$d$] (A) at (1,-0.35);
	\coordinate [label=center:$v$] (A) at (0,-0.35);
	%\draw (18:4) node{} ;
	
	\end{tikzpicture} \quad
	\begin{tikzpicture}[thick,scale=0.9]%
	
	%	\draw[dotted] (1.2,1.2)--(1,0)--(0.8,1.2);
	\draw[dashed] (-1,0) -- (0,0);
	\draw[dashed] (0,0) -- (1,0);
	\draw[dashed] (1,0)  -- (2,0);
	\draw node{} (0,0) node{} ;
	\draw node{} (1,0) node{} ;
	%\coordinate [label=center:$S$] (A) at (1,1.4);
	\coordinate [label=center:$d$] (A) at (1,-0.35);
	\coordinate [label=center:$v$] (A) at (0,-0.35);
	\end{tikzpicture} \quad
	\begin{tikzpicture}[thick,scale=0.9]%
	\draw[dotted] (0.2,1.2)--(0,0)--(-0.2,1.2);
	\draw[dashed] (-1,0) -- (0,0);
	\draw[dashed] (0,0) -- (1,0);
	\draw[dashed] (1,0)  -- (2,0);
	\draw node{} (0,0) node{} ;
	\draw node{} (1,0) node{} ;
	\coordinate [label=center:$S'$] (A) at (0,1.4);
	\coordinate [label=center:$d$] (A) at (1,-0.35);
	\coordinate [label=center:$v$] (A) at (0,-0.35);
	
	\end{tikzpicture}

	\caption{the graph $G$, $S$ being pruned from $G$ and $S$ being regrafted at $v$}
	\label{fig:SPR}
\end{figure}
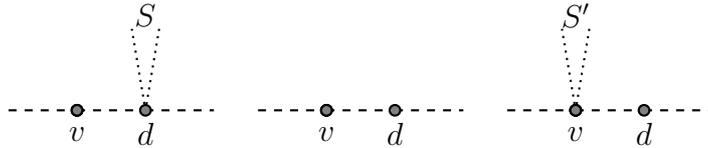

%%%%%%%%%%%%%%%%%%%%%%%%%%%%%%%%%%%%%%%
%%%%%%%%%%%%%%%%%%%%%%%%%%%%%%%%%%%%%%%
%%%%%%%% Preliminaries %%%%%%%%%%%%%%%%
%%%%%%%%%%%%%%%%%%%%%%%%%%%%%%%%%%%%%%%
%%%%%%%%%%%%%%%%%%%%%%%%%%%%%%%%%%%%%%%

%\section{Preliminaries}

%%%%%%%%%%%%%%%%%%%%%%%%%%%%%%%%%%%%%%%
%%%%%%%%%%%%%%%%%%%%%%%%%%%%%%%%%%%%%%%
%%%%%%%% Rect curv E4  %%%%%%%%%%%%%%%%
%%%%%%%%%%%%%%%%%%%%%%%%%%%%%%%%%%%%%%%
%%%%%%%%%%%%%%%%%%%%%%%%%%%%%%%%%%%%%%%

%\section{Rectifying curves in $\mathbb{E}^4$}

\section{Monotonicity in the matching number}\label{2}

In this section, we will prove the following proposition.

\begin{prop}\label{monotonicity}
	For fixed $n$, when $m_1<m_2$ and the sets $\U(n,m_1)$ and $\U(n,m_2)$ are both nonempty, then $W(\U(n,m_1)) < W(\U(n,m_2))$
\end{prop}

Assume this proposition is not true. In that case there exist some $n$ and $m$ such that 
$\U(n,m)$ and $\U(n,m+1)$ are both nonempty and $W(\U(n,m)) \ge  W(\U(n,m+1)).$
For some fixed $n$, we take the least integer $m$ for which this holds and take an extremal graph $G \in \U(n,m)$ with $W(G)=W(\U(n,m))$.
Note that $G$ cannot be a cycle itself, since then $\U(n,m+1)$ would be empty.
Hence $G$ contains some cycle $C_k$ where at each vertex $r_i$ of the cycle there are attached some trees $T_i$ (possibly consisting only of the vertex $r_i$).

Assume at least one of those trees attached is not a path. Then we can look to a vertex $w$ on such a tree with degree at least $3$ such that $d(w,C_k)$ is maximal.
There are at least two leaves such that the shortest paths from these leaves to $C_k$ contain $w$, by the choice of $w$.
We call two of them $l_1$ and $l_2$.

%%%Regraft?

Let $S$ be the path from $w$ to $l_1$. We prune $S$ and regraft it at $l_2$.

%We remove the path from $w$ to $l_2$ except from the vertex $w$ itself and connect the path to $l_1$.

The new graph will have matching number $m$ or $m+1$ and a Wiener index which is strictly larger than that of $G$.
This is a contradiction with $W(G)=W(\U(n,m)) \ge W(\U(n,m+1)).$

Hence all trees $T_i$ are paths.
Let $P_i=T_i \backslash \{r_i\}$.
If there is only one $P_i$ which is nonempty, then $G$ has a maximum matching, which cannot be the case since $\U(n,m+1)$ is nonempty.
Let $I$ be the set of indices $i$ such that $P_i$ is nonempty.
Take $i_1,i_2 \in I$ and wlog \begin{equation}\label{ass}
\sum_{t \in I \backslash \{i_1,i_2\}} d(r_{i_1},r_t) \lvert P_t \rvert \ge \sum_{t \in I \backslash \{i_1,i_2\}} d(r_{i_2},r_t) \lvert P_t \rvert
\end{equation}
Let $l_{i_1}$ be the leaf of $P_{i_1}$.
Note that $W(G)$ equals 
\begin{align}
\begin{split} \label{eq:calcW}
W(G)={}&\sum_{1\le t \le 2} \sum_{u,v \in P_{i_t}} d_G(u,v) +  \sum_{u \in P_{i_1}} \sum_{v \in P_{i_2}} d_G(u,v)+
\sum_{u,v \in G\backslash \{P_{i_1},P_{i_2}\}} d_G(u,v)  \\
&  +\sum_{u \in  \{P_{i_1},P_{i_2}\} } \sum_{v \in C_k } d_G(u,v)+\sum_{u \in  \{P_{i_1},P_{i_2}\} } \sum_{t \in I \backslash \{i_1,i_2\} } \sum_{v \in P_{i_t}} d_G(u,v)  \end{split}
\end{align}

Prune $P_{i_2}$ and regraft it at $l_{i_1}$, to get a graph $G'$. This increases the matching number at most by $1$.
Using \eqref{eq:calcW}, we get that 
\begin{align*}
W(G')-W(G)= &0  - \left( d(r_{i_1},r_{i_2}) +1\right) \lvert P_{i_1} \rvert \lvert P_{i_2} \rvert +0 +\lvert G \backslash \{P_{i_1},P_{i_2}\} \rvert \lvert P_{i_1} \rvert \lvert P_{i_2} \rvert + \\
& \sum_{t \in I \backslash \{i_1,i_2\}} \left( d(r_{i_1},r_t) - d(r_{i_2},r_t) \right) \lvert P_t \rvert  \lvert P_{i_2} \rvert \\
&\ge \left( k- d(r_{i_1},r_{i_2}) -1 \right)  \lvert P_{i_1} \rvert \lvert P_{i_2} \rvert \\
&> 0
\end{align*}
%elementen van P_{i_2} worden een afstand $|P_{i_1}|$ verder van C_k geplaatst en na verplaatsing van r_2 naar r_1 ook over die afstand van de andere P_t
since $d(r_{i_1},r_{i_2}) \le k-2$, $ \lvert G \backslash \{P_{i_1},P_{i_2}\} \rvert \ge \lvert C_k \rvert =k$ and by \eqref{ass}.
Hence $G'$ satisfies $W(G')>W(G)$, while the matching number of $G'$ is at most $m+1$. Since the assumption was that $G$ had the largest Wiener index among all graphs on $n$ vertices with matching number no larger than $m+1$, we get a contradiction.
Hence the assumption was false and no such graph $G$ did exist, which proves Proposition~\ref{monotonicity}. %\qedsymbol 

\section{Reducing the cycle length}\label{3}

In this section, we prove that the extremal graphs cannot contain a cycle of length at least $5$, as stated in the following proposition.

\begin{prop}\label{C3C4}
	Given a graph in $\U(n,m)$ with Wiener index $W(\U(n,m))$, the cycle of $G$ will be $C_3$ or $C_4$.
\end{prop}

\begin{figure}[h]
	\centering
	\begin{tikzpicture}[thick,scale=0.66]%
	\draw \foreach \x in {90,162,234,306} {
		(\x:4) node{} -- (\x+72:4)

	};
	\draw (18:4) node{} ;
	
	\coordinate [label=left:$T_{k-1}$] (A) at (18:5.4); 
	\coordinate [label=left:$T_3$] (A) at (90:4.25); 
	\coordinate [label=left:$T_2$] (A) at (162:4.25); 
	\coordinate [label=left:$T_1$] (A) at (234:4.25); 
	\coordinate [label=left:$T_k$] (A) at (306:4.5); 
	
	\coordinate [label=left:$r_{k-1}$] (A) at (18:3.65); 
	\coordinate [label=left:$r_3$] (A) at (85:3.5); 
	\coordinate [label=left:$r_2$] (A) at (162:3.25); 
	\coordinate [label=left:$r_1$] (A) at (240:3.25); 
	\coordinate [label=left:$r_k$] (A) at (310:3.5); 
	
	%\draw (18:4.5) {$T_1$};
	\draw[dotted] (18:4) 
	arc (18:90:4);	
	
	\end{tikzpicture}
	\caption{graph $G$ containing $C_k$}
	\label{fig:C_k graph}
\end{figure}
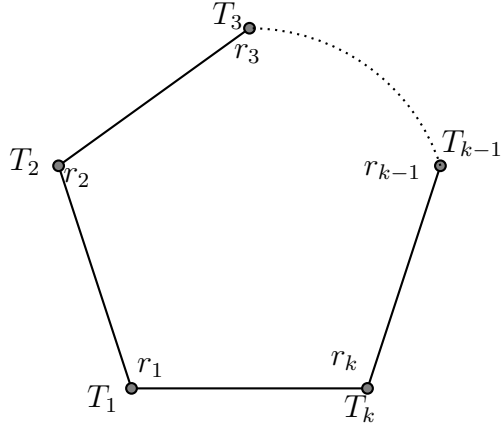

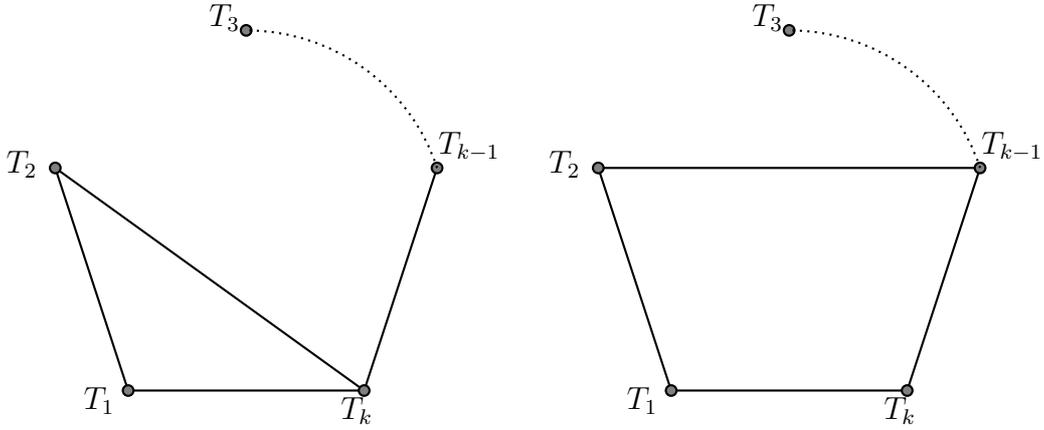
\begin{figure}[h]
	\centering
	
	\begin{tikzpicture}[thick,scale=0.66]%
	\draw \foreach \x in {162,234,306} {
		(\x:4)  -- (\x+72:4)

	};
	\draw (162:4)--(306:4);
	\draw \foreach \x in {18,90,162,234,306} {
		(\x:4) node{}
		
	};
	
	\coordinate [label=left:$T_{k-1}$] (A) at (18:5.4); 
	\coordinate [label=left:$T_3$] (B) at (90:4.25); 
	\coordinate [label=left:$T_2$] (C) at (162:4.25); 
	\coordinate [label=left:$T_1$] (D) at (234:4.25); 
	\coordinate [label=left:$T_k$] (E) at (306:4.5); 
	
	%\draw (18:4.5) {$T_1$};
	\draw[dotted] (18:4) 
	arc (18:90:4);	
	
	\end{tikzpicture} \quad
	\begin{tikzpicture}[thick,scale=0.66]%
	\draw \foreach \x in {162,234,306} {
		(\x:4)  -- (\x+72:4)

	};
	\draw (162:4)--(18:4);
	\draw \foreach \x in {18,90,162,234,306} {
		(\x:4) node{}
		
	};
	
	\coordinate [label=left:$T_{k-1}$] (A) at (18:5.4); 
	\coordinate [label=left:$T_3$] (A) at (90:4.25); 
	\coordinate [label=left:$T_2$] (A) at (162:4.25); 
	\coordinate [label=left:$T_1$] (A) at (234:4.25); 
	\coordinate [label=left:$T_k$] (A) at (306:4.5); 
	
	%\draw (18:4.5) {$T_1$};
	\draw[dotted] (18:4) 
	arc (18:90:4);	
	
	\end{tikzpicture}

	\caption{Modifications $G_1$ and $G_2$ for graph $G$}
	\label{fig:C_k graph1}
\end{figure}

Assume some extremal graph $G \in \mathbb{U}(n,m)$ contains a cycle $C_k$ of length $k\ge 5$.
Then there is some tree $T_i$ attached to every vertex $r_i$ of $C_k$.
We will call these rooted trees $T_1$ up to $T_k$ and their roots which lie on the $k$-gon $C_k$ will be called $r_1$ up to $r_k$ in cycle order. Note that it is possible that such a tree equals a single vertex, the root.
With $|T_i|$ we will denote the number of vertices of the tree $T_i$ and we assume
\begin{equation}\label{T3max}
|T_3|=\max\{|T_1|,|T_2|,\ldots,|T_k|\}
\end{equation}
In Figure~\ref{fig:C_k graph} we see the graph $G$ consisting of a $k$-cycle with the trees $T_i$ attached to the vertices of the $k$-cycle.
In Figure~\ref{fig:C_k graph1} there are drawn two possible modifications $G_1, G_2$ for this graph.

In words, $G_1$ (resp. $G_2$) is obtained from $G$ by deleting the edge $r_2r_3$ and adding $r_2r_k$ (resp. $r_2r_{k-1}$).
Clearly, $G_1$ and $G_2$ gave the same number of edges as $G$ does.

We will prove that at least one of $G_1, G_2$ has a higher Wiener index and a matching number which is at most the matching number of $G$, implying that $G$ was not an extremal graph.

First we derive that both modifications have a higher Wiener index.
Note that $$W(G)=\sum_{1\le i \le k} \sum_{u,v \in T_i} d_G(u,v) + \sum_{1 \le i<j \le k} \sum_{u \in T_i} \sum_{v \in T_j} d_G(u,v)$$
which implies that 
\begin{equation}\label{WG2-WG} W(G_2)-W(G)=\sum_{1 \le i<j \le k} \lvert T_i \rvert  \lvert T_j \rvert \left( d_{G_2}(r_i,r_j)-d_G(r_i,r_j) \right)\end{equation}
and 
\begin{equation}\label{WG1-WG2}
W(G_1)-W(G_2)=\sum_{1 \le i<j \le k} \lvert T_i \rvert  \lvert T_j \rvert \left( d_{G_1}(r_i,r_j)-d_{G_2}(r_i,r_j) \right).
\end{equation}
We will prove that $W(G)<W(G_2) < W(G_1)$.

Note that  $d_{G_2}(r_i,r_j) \le d_{G_1}(r_i,r_j)$ for all $1 \le i < j \le k$ except for $i=2,j=k$.
By \eqref{WG1-WG2}, we get that $W(G_1)-W(G_2)\ge \lvert T_2 \rvert  \lvert T_{k-1} \rvert + \lvert T_2 \rvert  \lvert T_3 \rvert -\lvert T_2 \rvert  \lvert T_k \rvert \ge 1$ since $\lvert T_i \rvert \ge 1$ for all $1\le i \le k$ and due to \eqref{T3max}.

Note that $d_{G_2}(r_i,r_j) \ge d_G(r_i,r_j)$ when $i$ and $j$ are both different from $2$, since there exists a shortest path in $G_2$ from $r_i$ to $r_j$ which is a subpath of the path $r_1r_kr_{k-1}\ldots r_3$ which also exists in $G$.

Note that 
\begin{equation}\label{pos}
d_{G_2}(r_i,r_3)=d_{G}(r_i,r_3)+(k-4) \mbox{ for }i \in \{1,2\}.
\end{equation}

Next, we have to do a small case distinction between $k$ even and $k$ odd.

When $k$ is odd, then $d_{G_2}(r_2,r_j)=d_G(r_2,r_j)-2$ for $\frac{k+5}{2}\le j\le k-1$ and $d_{G_2}(r_2,r_j)=d_G(r_2,r_j)-1$ for $j=\frac{k+3}{2}$.
For the remaining values of $j$, one has $d_{G_2}(r_2,r_j) \ge d_G(r_2,r_j)$.
Using \eqref{WG2-WG} and \eqref{pos}, we have for $k$ odd that 
\begin{align*}
W(G_2)-W(G) &\ge (k-4) \left( \lvert T_1 \rvert + \lvert T_2 \rvert\right))\lvert T_3 \rvert 
-2 \sum_{\frac{k+5}{2} \le j \le k-1} \lvert T_2 \rvert  \lvert T_j \rvert 
- \lvert T_2 \rvert  \lvert T_{\frac{k+3}{2}} \rvert\\
&\ge (k-4)  \lvert T_1 \rvert  \lvert T_3 \rvert >0
\end{align*} due to \eqref{T3max}, the fact that $\lvert T_1 \rvert >0$ and the fact that there are $\frac{k-5}{2}$ integral numbers in the interval $[\frac{k+5}{2},k-1]$.
When $k$ is even, $d_{G_2}(r_2,r_j)=d_G(r_2,r_j)-2$ for $2+\frac{k}{2}\le j\le k-1$, and $d_{G_2}(r_2,r_j) \ge d_G(r_2,r_j)$ for the remaining values of $j$.
Using \eqref{WG2-WG} and \eqref{pos}, we have for $k$ even that 
\begin{align*}
W(G_2)-W(G) &\ge (k-4) \left( \lvert T_1 \rvert + \lvert T_2 \rvert\right) \lvert T_3 \rvert 
-2 \sum_{\frac{k+4}{2} \le j \le k-1} \lvert T_2 \rvert  \lvert T_j \rvert 
&\ge (k-4)  \lvert T_1 \rvert  \lvert T_3 \rvert >0
\end{align*} by \eqref{T3max}, $\lvert T_1 \rvert >0$ and the fact that the sum is over $\frac{k-4}{2}$ integers.

Next, we show that the matching number of $G$ is not smaller than the matching numbers of both $G_1$ and $G_2$.
Assume to the contrary that the matching number of $G$ is strictly smaller than the matching numbers of $G_1$ or $G_2$.
It is clear that a maximal matching in $G$ can be modified to a matching such that the submatchings  are optimal for every tree $T_i$ and do not use $r_i$ if there exists a maximal matching for $T_i$ without using $r_i$.
So starting from the optimal matchings in every $T_i$ such that $r_i$ is not used when not necessary, the remaining task is to find an optimal matching between the $r_i$ which are not used.
Since the optimal matching in $G_1$ is strictly larger than the optimal matching in $G$, we will use $r_2 r_k$ in that matching and so $r_1$ has to be used in the optimal matching of $T_1$, since otherwise we could use $r_2 r_1$ instead of $r_2 r_k$ in $G$ and take the other pairs as in the matching of $G_1$.
This implies also that $r_2$ and $r_k$ are not used in the optimal matchings of $T_2$ and $T_k$ respectively.
Similarly we will use $r_2 r_{k-1}$ in the optimal covering of $G_2$.
But since $r_1$ is used in the optimal covering of $T_1$ and $r_k$ is not used in the optimal covering of $T_k$, we can replace $r_2 r_{k-1}$ by $r_k r_{k-1}$ and so the same matching would work for $G$, contradiction.

We conclude that $G_1$ or $G_2$ have matching numbers not larger than that of $G$, but have a Wiener index which is greater than the one of $G$.
Together with Proposition~\ref{monotonicity}, this implies that it is impossible that $W(G)=W(\U(n,m))$. So the assumption at the beginning was wrong, implying that no extremal graph can contain a cycle of length at least $5$, proving the proposition.

%$$W(G_m)-W(G)=\sum_{1 \le i<j \le k} \lvert T_i \rvert  \lvert T_j \rvert \left( d_{G_m}(r_i,r_j)-d_G(r_i,r_j) \right).$$
%Next, we note that $d_{G_m}(r_i,r_j) \ge d_G(r_i,r_j)$ when $i$ and $j$ are both different from $2$.

\section{Optimal form of trees}\label{4}

In this section, we will prove that each tree $T_i$ of an extremal graph is a concatenation of a path and a star, see %i.e. the optimal configurations will be as shown in 
figure~\ref{fig:pos_con}.

\begin{prop}
	For a graph $G$ with $W(G)=W(\U(n,m))$, $G$ is a cycle $C_3$ or $C_4$ with trees attached to it, each of which is a path attached to a star.
\end{prop}

We know already by Proposition~\ref{C3C4} that an extremal graph $G$ is a cycle $C_3$ or $C_4$ with some trees $T_i$ attached.

%We may assume wlog that $\lvert T_1 \rvert \ge \lvert T_i \rvert $ for every $i$.

For every $i$, we take a longest path $\P$ in $T_i$ starting from $r_i$ and call a leaf of that longest path $l_i$ and the adjacent vertex to $l_i$ on that path $c_i$.

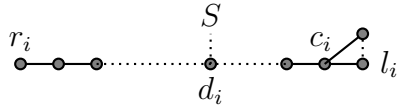
\begin{figure}[h]
	\centering
	
	\begin{tikzpicture}[thick,scale=1]%
	
	\draw \foreach \x in {0,0.5,3.5,4} {
		%\draw[dotted] (\x:2.75)--(\x:3.5)  
		(0:\x)  node {} --(0:\x+0.5)
	};
	
	\draw[dotted] (4.5,0)--(4.5,0.4);
	\draw (0:4.5)  node {};
	\draw (0:4)   --(4.5,0.4) node {};
	
	\draw (0:1)  node {};
	\draw (0:2.5)  node {};

	\draw[dotted] 	 (0:1)--(0:3.5);
	\draw[dotted]  (2.5,0)--(2.5,0.4);

	\coordinate [label=center:$l_i$] (A) at (0:4.85); 
	\coordinate [label=center:$r_i$] (A) at (0,0.3); 
	\coordinate [label=center:$S$] (A) at (2.5,0.65); 
	\coordinate [label=center:$c_i$] (A) at (3.95,0.3);
	\coordinate [label=center:$d_i$] (A) at (2.5,-0.35);
	
	\end{tikzpicture} 
	%\quad
	%	\begin{tikzpicture}[thick,scale=0.8]%
	%
	%		\draw \foreach \x in {0,0.5,3.5,4} {
	%		%\draw[dotted] (\x:2.75)--(\x:3.5)  
	%		(90:\x)  node {} --(90:\x+0.5)
	%	};
	%	
	%	\draw[dotted] (0,4.5)--(0.8,4.5);
	%	\draw (90:4.5)  node {};
	%	\draw (90:4)   --(0.8,4.5) node {};
	%	
	%	\draw (90:1)  node {};
	%	\draw (90:2.5)  node {};
	%	
	%	
	%	\draw[dotted] 	 (90:1)--(90:3.5);
	%	\draw[dotted]  (0,2.5)--(0.4,2.5);
	%	
	%	
	%	
	%	
	%	
	%	\coordinate [label=center:$l_1$] (A) at (90:4.85); 
	%	\coordinate [label=center:$l_a$] (A) at (0.8,4.85);
	%	\coordinate [label=center:$r_1$] (A) at (0.3,0); 
	%	\coordinate [label=right:$S$] (A) at (80:2.5); 
	%	\coordinate [label=center:$A_j$] (A) at (-0.375,4);
	%	\coordinate [label=center:$A_p$] (A) at (-0.375,2.5);
	%	\coordinate [label=center:$A_2$] (A) at (-0.375,1);
	%	\coordinate [label=center:$A_1$] (A) at (-0.375,0.5);
	%	
	%	\end{tikzpicture}
	%	

	\caption{subtree $T_i$}
	\label{fig:paths}
\end{figure}

Assume that some tree $T_i$ is not a path attached to a star. 
Then there exists a nonempty rooted subtree $S$ with root $d_i$ (possibly equal to $r_i$) on $\P$ which is closest to $c_i$, as shown in Figure~\ref{fig:paths}.
Here $\P \cap S= d_i$ and $\P \cup S$ contain all edges adjacent to $d_i$.
%Removing $S$ from its initial connection point $d_i$ to the point $c_i$ will strictly increase the Wiener index of the graph as $\lvert T_1 \rvert \ge \lvert T_i \rvert $. On the other hand the matching number cannot increase by this replacement of $S$ and so .
%So this implies the trees $T_i$ with $i \ge 2$ have to be of the form we mention.

The vertex $d_i$ partitions the edge set of $G$ into three parts: $E(H_1)$, $E(H_2)$ and $E(S)$ where $H_1$ is a tree containing $l_i$ and $d_i$ as leaves, $H_2$ is a unicyclic subgraph with $d_i$ being a leaf.

Let $V_1=V(H_1)\backslash \{d_i\}$ and $V_2=V(H_2)\backslash \{d_i\}$ and $V_S= V(S)\backslash \{d_i\}$.

If $\lvert V_1 \rvert \le \lvert V_2 \rvert$, 
then pruning $S$ from $G$ and regrafting $S$ at $c_i$ will strictly increase the Wiener index.
Note for this that $\sum_{u \in V_S,v \in V_1} d(u,v)$ has decreased with at most $d(d_i,c_i) \lvert V_1 \rvert \lvert V_S \rvert $ and 
$\sum_{u \in V_S,v \in V_2} d(u,v)$ has increased with $d(d_i,c_i) \lvert V_2 \rvert  \lvert V_S \rvert $, while $\sum_{u \in V_S} d(u,d_i)$ strictly increases by the value $\lvert V_S \rvert d(d_i,c_i)$.

Also, the matching number has not been increased.
Doing this, we see the trees attached to the cycle, with possible one exception, are each a path attached to a star.

If  $\lvert V_1 \rvert > \lvert V_2 \rvert$, then analogously pruning $S$ and regrafting it at $c_j$ (with $j\not=i$) will do the job (assuming the subtree $T_j$ has length at least $1$).

In the case there was only one tree attached to the cycle, we cannot do this and we have to use other replacements.
If the cycle is $C_4$, we can prune $S$ and regraft it at  $r_j$ (with $r_j$ and $r_i$ being opposite corners of $C_4$).

So from now on, we assume the cycle is $C_3.$

In the case there are multiple subtrees attached to $\P$ with the root not equal to $c_i$, if there is any subtree $S'$ such that there is some maximum matching that does not use its root $d'_{i} \in \P$, we can prune $S'$ and regraft is at some $r_j$ ($j \not = i$) and conclude.

we will assume $S'$ is the rooted subtree with root $d'_{i}$ which is second closest to $c_i.$
We can prune $S$ and regraft it at $d'_{i}$, so the Wiener index increases, while the matching number does not and conclude again that the original graph was not extremal.

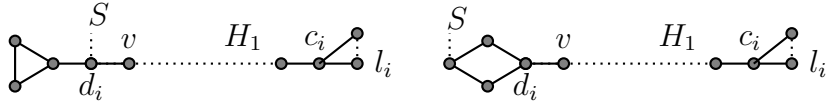
\begin{figure}[h]
	\centering
	
	\begin{tikzpicture}[thick,scale=1]%
	
	\draw[dotted] 	 (0:1)--(0:3.5);
	\draw[dotted]  (1,0)--(1,0.4);
	
	\draw (0.5,0)-- (0,-0.3)--(0,0.3)--(0.5,0);
	\draw \foreach \x in {0.5,1,3.5,4} {
		%\draw[dotted] (\x:2.75)--(\x:3.5)  
		(0:\x)  node {} --(0:\x+0.5)
	};

	\draw[dotted] (4.5,0)--(4.5,0.4);
	\draw (0:4.5)  node {};
	\draw (0:4)   --(4.5,0.4) node {};
	
	\draw (0:1)  node {};
	\draw (0:1.5)  node {};
	%	\draw (0:2.5)  node {};

	\draw (0,0.3) node{};
	\draw (0,-0.3) node{};

	\coordinate [label=center:$d_i$] (A) at (1,-0.3); 
	\coordinate [label=center:$l_i$] (A) at (0:4.85); 
	\coordinate [label=center:$v$] (A) at (1.5,0.3); 
	\coordinate [label=center:$S$] (A) at (1.1,0.65); 
	\coordinate [label=center:$c_i$] (A) at (3.95,0.3);
	\coordinate [label=center:$H_1$] (A) at (3,0.35);
	
	\end{tikzpicture} 
	\quad
	\begin{tikzpicture}[thick,scale=1]%

	\draw[dotted] 	 (0:1)--(0:3.5);
	\draw[dotted]  (0,0)--(0,0.4);
	
	\draw (0,0)-- (0.5,-0.3)--(1,0)--(0.5,0.3)--(0,0);
	\draw \foreach \x in {1,3.5,4} {
		%\draw[dotted] (\x:2.75)--(\x:3.5)  
		(0:\x)  node {} --(0:\x+0.5)
	};

	\draw[dotted] (4.5,0)--(4.5,0.4);
	\draw (0:4.5)  node {};
	\draw (0:4)   --(4.5,0.4) node {};
	
	\draw (0:1)  node {};
	\draw (0:1.5)  node {};
	\draw (0,0) node{};
	%	\draw (0:2.5)  node {};

	\draw (0.5,0.3) node{};
	\draw (0.5,-0.3) node{};

	\coordinate [label=center:$l_i$] (A) at (0:4.85); 
	\coordinate [label=center:$d_i$] (A) at (1,-0.3); 
	\coordinate [label=center:$S$] (A) at (0.1,0.6); 
	\coordinate [label=center:$c_i$] (A) at (3.95,0.3);
	\coordinate [label=center:$v$] (A) at (1.5,0.3); 
	\coordinate [label=center:$H_1$] (A) at (3,0.35);
	\end{tikzpicture}

	\caption{Constructing a better graph}
	\label{fig:change_configuration}
\end{figure}

In the final case, we assume $S$ is the only such subtree connected to $\P$.
If $d_i=r_i$, we can prune and regraft $S$ at $r_j$.
If $d(d_i,r_i)>2$, we prune and regraft it at a vertex $v$ which is an even distance closer to $r_i$, which also increases the Wiener index while the matching number is constant. 
If $d(d_i,r_i)=1$, we easily can compare the configuration with an other one using a $C_4$ as cycle, as shown in Figure~\ref{fig:change_configuration}.
Here $v$ is the neighbour of $d_i$ which is closer to $c_i.$
The matching number of both configurations is the same, while the Wiener index increases with $\left( 2\lvert V_S \rvert -1 \right) \lvert V_1 \rvert -\lvert V_S \rvert>0,$ since $\lvert V_S \rvert \ge 1$ and $\lvert V_1 \rvert >3$.

Since removing $S$ has increased the Wiener index, while the matching number has not been increased, we have a contradiction in our assumption that $W(G)=W(\U(n,m))$ due to Proposition~\ref{monotonicity}.

\section{Only one long tree}\label{5}

\begin{figure}[h]
	\centering
	
	\begin{tikzpicture}[thick,scale=0.8]%
	\draw \foreach \x in {0,120,240} {
		%(\x:3) node {} -- (\x+120,3))
		(\x:2) node {}  -- (\x+120:2)
		(\x:2)--(\x:2.75)  node {}
	};
	\draw[dotted] 	\foreach \x in {0,120,240} {
		(\x:2.75)--(\x:3.75) 
	};
	
	\draw \foreach \a in {6,-6} \foreach \x in {0,120,240} {
		(\x:3.75)--(\x+\a:4.25) node {}
	};
	\draw[dotted]  \foreach \x in {0,120,240} {
		(\x-6:4.25)--(\x+6:4.25) 
	};
	
	\draw \foreach \x in {0,120,240} {
		%\draw[dotted] (\x:2.75)--(\x:3.5)  
		(\x:3.75)  node {}
		
	};

	\coordinate [label=left:$A_0$] (A) at (10.1:2.3); 
	\coordinate [label=right:$B_0$] (A) at (117.5:2.15); 
	\coordinate [label=right:$C_0$] (A) at (242.5:2.15); 
	\coordinate [label=left:$A_1$] (A) at (6.75:3); 
	\coordinate [label=right:$B_1$] (A) at (117.5:2.8); 
	\coordinate [label=right:$C_1$] (A) at (242.5:2.8); 
	\coordinate [label=left:$A_j$] (A) at (5.5:3.95); 
	\coordinate [label=right:$B_k$] (A) at (117.5:3.7); 
	\coordinate [label=right:$C_l$] (A) at (242.5:3.7);

	\coordinate [label=center:$a$] (A) at (0:4.5); 
	\coordinate [label=center:$b$] (A) at (120:4.5); 
	\coordinate [label=center:$c$] (A) at (240:4.5); 
	
	\end{tikzpicture} \quad
	\begin{tikzpicture}[thick,scale=0.8]%
	\draw \foreach \x in {0,90,180,270} {
		%(\x:3) node {} -- (\x+90,3))
		(\x:2) node {}  -- (\x+90:2)
		(\x:2)--  (\x:2.75)  node {}
	};
	\draw[dotted] \foreach \x in {0,90,180,270} {
		(\x:2.75)--(\x:3.75)  
	};

	\draw \foreach \a in {5,-5} \foreach \x in {0,90,180,270} {
		(\x:3.75)--(\x+\a:4.25) node {}
	};
	\draw \foreach \x in {0,90,180,270} {
		%\draw[dotted] (\x:2.75)--(\x:3.5)  
		(\x:3.75)  node {}
		
	};
	
	\draw[dotted]  \foreach \x in {0,90,180,270} {
		(\x-5:4.25)--(\x+5:4.25) 
	};
	
	\coordinate [label=left:$A_0$] (A) at (2.5,0.35); 
	\coordinate [label=left:$B_0$] (A) at (0.85,2.1); 
	\coordinate [label=left:$C_0$] (A) at (-1.8,0.35); 
	\coordinate [label=left:$D_0$] (A) at (0.8,-2.1); 
	\coordinate [label=left:$A_1$] (A) at (3.15,0.35); 
	\coordinate [label=left:$B_1$] (A) at (0.85,2.8); 
	\coordinate [label=left:$C_1$] (A) at (-2.55,0.35); 
	\coordinate [label=left:$D_1$] (A) at (0.85,-2.8); 	
	\coordinate [label=left:$A_j$] (A) at (3.85,0.3); 
	\coordinate [label=left:$B_k$] (A) at (0.85,3.7); 
	\coordinate [label=left:$C_l$] (A) at (-3.4,0.35); 
	\coordinate [label=left:$D_h$] (A) at (0.85,-3.7); 
	
	\coordinate [label=center:$a$] (A) at (0:4.5); 
	\coordinate [label=center:$b$] (A) at (90:4.5); 
	\coordinate [label=center:$c$] (A) at (180:4.5); 
	\coordinate [label=center:$d$] (A) at (270:4.5); 
	
	\end{tikzpicture}

	\caption{possible configurations}
	\label{fig:pos_con}
\end{figure}
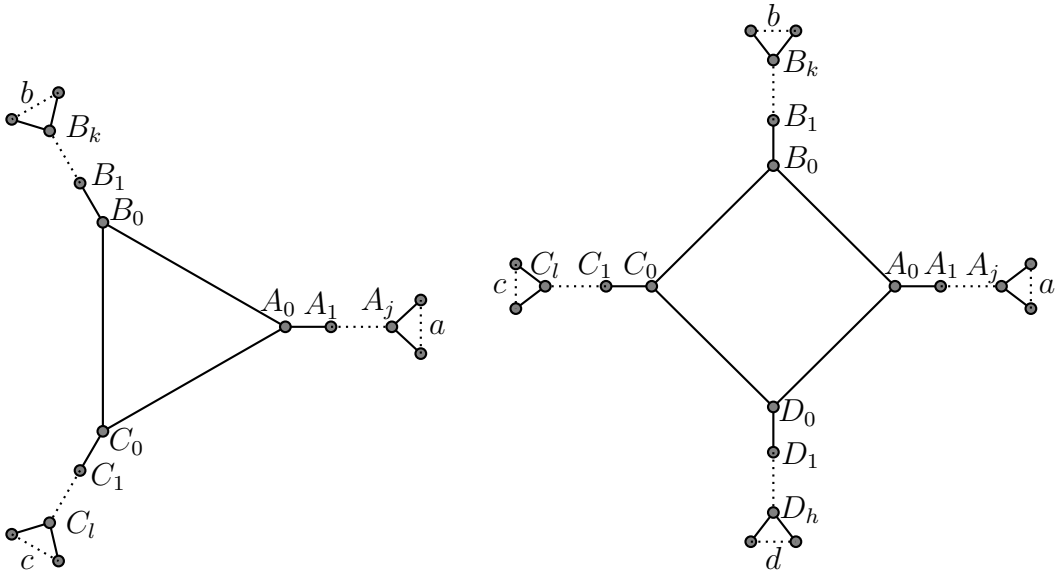

From the previous sections, we can conclude that the possible extremal graphs are as represented in Figure~\ref{fig:pos_con}, here $a,b,c,d$ are the numbers of leaves in every subtree.
In this section, using some calculations we prove some extra conditions on the configuration of an extremal graph.

\begin{prop}
	An extremal graph in $\U(n,m)$ has at most one tree with height larger than $1$ connected to the cycle.
	When the cycle is $C_3$, then the extremal graph is isomorphic to a graph $G^3_{a,b,c,j}$  (see Figure~\ref{fig:G3abcj_G4acj}) with $a \ge b \ge c$.
	Furthermore, when the cycle is $C_4$, the extremal graph is isomorphic to a graph $G^4_{a,c,j}$ (see Figure~\ref{fig:G3abcj_G4acj}) where $a \in \{c,c+1\}$.
	
\end{prop}

\subsection{Case $G$ contains $C_3$}

Using the fact that $W(P_q)=\binom{q+1}{3}$, a calculation gives that the Wiener index of the first configuration in Figure~\ref{fig:pos_con} equals
\begin{align*}
W(G)=&\binom{k+l+3}{3}+\binom{j+k+3}{3}+\binom{l+j+3}{3}-\binom{k+2}{3}-\binom{l+2}{3}-\binom{j+2}{3}\\
&+2\left( \binom{a}{2}+\binom{b}{2}+\binom c2 \right)
+ac(l+j+3)+ab(k+j+3)+bc(k+l+3)\\
&+(a+b+c)\left( \binom{j+2}{2}+\binom{k+2}{2}+\binom{l+2}{2}\right)
\\&+a(k+l+2)(j+1)+b(k+1)(l+j+2)+c(l+1)(k+j+2)
\end{align*}
Wlog $a=\max\{a,b,c\}$.
We will prove that $W(G')>W(G)$ when $G'$ is the graph determined by $j'=j+k+l$ and $k'=l'=0$ and $(a',b',c')=(a,b,c)$, if at most one of $j,k,l$ equals zero. 
This is a consequence of the following calculation.

\begin{align*}
W(G')-W(G)%=&jkl+jk+jl+kl\\
%&+k(a-b)c+l(a-c)b\\
%&+(a+b+c)(jk+jl+kl)\\
%&-a(kj+lj-k-l)-b(kl+kj+k)-c(kl+lj+l)\\
=&jkl+jk+jl+kl+k(a-b)c+l(a-c)b\\
&+akl+bjl+cjk+k(a-b)+l(a-c)\\
>& 0
\end{align*}
since $a-b,a-c\ge 0$ and at least one of $jk,jl,kl$ is strictly positive.

We see that the matching numbers of $G'$ and $G$ are the same.
Note we can only use one leaf of the star at the end.
Removing $\max\{a-1,0\}$ edges from the right star and similarly for the two other stars, we have a triangle with some attached paths.
Now a maximum matching uses all vertices, or all vertices minus one when the total number is odd.
The same holds after the operation $j+k+l\rightarrow j'$.
This implies that the extremal graphs containing $C_3$ have two trees of height at most $1$.
Due to the terms $k(a-b)$ and $l(a-c)$, the new configuration is strictly better if $a>b$ and $k>0$ or $a>c$ and $l>0$ originally.
Hence $a= \max \{a,b,c\}$ is strictly necessary when we take $j>0$ in the graph $G^3_{a,b,c,j}.$

\subsection{Case $G$ contains $C_4$}\label{sec:C4}

We calculate the Wiener index for the graph in function of the parameters $a,b,c,d,h,j,k,l$.
We use formulas like $W(P_q)=\binom{q+1}{3}$,
\begin{align*}
\sum_{0\le x \le j, 0\le y \le k} d(A_x,B_y) &=  \binom{k+j+3}{3}-\binom{j+2}{3}-\binom{k+2}{3}\\&= %(k+1)\binom{j+1}{2}+(j+1)\binom{k+1}{2}+(k+1)(j+1) \\&=
(k+1)(j+1)\frac{k+j+2}{2}
\end{align*} 
and
\begin{align*}
\sum_{0\le x \le h, 0\le y \le k} d(D_x,B_y) &=  %\binom{k+j+3}{3}-\binom{j+2}{3}-\binom{k+2}{3}\\&=
(k+1)\binom{h+1}{2}+(h+1)\binom{k+1}{2}+2(k+1)(h+1) \\&=
(k+1)(h+1)\frac{k+h+4}{2}
\end{align*}
multiple times.

\begin{align*}
W(G)=&\binom{h+2}{3}+\binom{j+2}{3}+\binom{k+2}{3}+\binom{l+2}{3}\\
&+(k+1)(j+1)\frac{k+j+2}{2}+(k+1)(l+1)\frac{k+l+2}{2}+(k+1)(h+1)\frac{k+h+4}{2}\\&+(l+1)(j+1)\frac{l+j+4}{2}+(h+1)(j+1)\frac{h+j+2}{2}+(h+1)(l+1)\frac{l+h+2}{2}\\
&+ab(k+j+3)+bc(k+l+3)+cd(l+h+3)+ad(h+j+3)\\&+ac(j+l+4)+bd(k+h+4)\\
&+2\left( \binom{a}{2}+\binom{b}{2}+\binom c2+\binom d2 \right)   \\
&+(a+b+c+d)\left(   \binom{h+2}{2}+\binom{j+2}{2}+\binom{k+2}{2}+\binom{l+2}{2} \right)\\
&+a(j+1)(k+h+2)+b(k+1)(l+j+2)+c(l+1)(k+h+2)+d(h+1)(l+j+2)\\
&+a(j+2)(l+1)+b(k+2)(h+1)+c(l+2)(j+1)+d(h+2)(k+1).\end{align*}

Given a graph $G$ with parameters $\{a,b,c,d,h,j,k,l\},$ we construct the graph $G'$ with parameters $\{\frac{a+b+c+d+\epsilon}{2},0, \frac{a+b+c+d-\epsilon}{2},0,0,h+j+k+l,0,0\},$
where $\epsilon \in \{0,1\}$ is chosen such that $\epsilon \equiv a+b+c+d \pmod 2$.
We will prove that $m(G') \le m(G)$ and $W(G') \ge W(G)$.

Some elementary arithmetic operations show the following (which can be checked for instance by a standard symbolic manipulation program):

\begin{align*}
W(G')-W(G)=
&hjk+hjl+hkl+jkl+2hj+hk+2hl+2jk+jl+2kl\\
&+ahk+ahl+akl+bhj+bhl+bjl+chj+chk+cjk+djk+djl+dkl\\&+ah+ak+bj+bl+ch+ck+dj+dl\\
&+\frac14 \left((d-a-b-c-1)^2-(\epsilon-1)^2)\right)h\\&+\frac14 \left((a-b-c-d-1)^2-(\epsilon-1)^2)\right)j\\
&+\frac14\left((b-a-c-d-1)^2-(\epsilon-1)^2)\right)k\\
&+\frac14\left((c-a-b-d-1)^2-(\epsilon-1)^2)\right)l\\
&+\frac{1}{2} \left((a-c)^2+(b-d)^2-\epsilon^2\right)\\&\ge0
\end{align*}
The inequality holds since every term is positive, i.e. when $\epsilon=1$, then $\max\{\lvert a-c \rvert, \lvert b-d\rvert\} \ge 1$ since $a+b+c+d$ is odd and when $\epsilon=0$, then $\lvert b-a-c-d-1 \rvert \ge 1$ since $a+b+c+d+1$ is odd and similarly for the three other differences of squares.
Equality holds if and only if every term is equal to zero.
This implies that at least three values in $\{h,j,k,l\}$ are zero.
If $h=j=k=l=0$, equality holds iff $\lvert a-c \rvert + \lvert b-d \rvert \le 1$.
When at least three values of $a,b,c,d$ are nonzero, $m(G)\ge3$ while $m(G')=2$ and so $G$ was not optimal by Proposition~\ref{monotonicity}.
In the other case, we note that $G \sim G'$.
In the case one value in $\{h,j,k,l\}$ is nonzero, wlog $j>0$, we get $b=d=0$, $\lvert a-c-1 \rvert = \lvert \epsilon -1 \rvert$ and $\lvert a-c \rvert = \lvert \epsilon \rvert$ implying $G' \sim G$ again.

\begin{figure}[h]
	\centering
	
	\begin{tikzpicture}[thick,scale=1]%
	\draw \foreach \x in {0} {
		%(\x:3) node {} -- (\x+120,3))
		(\x:2) node {}  -- (\x+120:2)
		(\x:2)--(\x:2.75)  node {}
	};
	\draw (120:2) node {}  -- (240:2);
	\draw (240:2) node {}  -- (0:2);
	
	\draw[dotted] 	\foreach \x in {0} {
		(\x:2.75)--(\x:3.75) 
	};
	
	\draw \foreach \a in {8,-8} \foreach \x in {120,240} {
		(\x:2)--(\x+\a:2.75) node {}
	};
	\draw[dotted]  \foreach \x in {120,240} {
		(\x-8:2.75)--(\x+8:2.75) 
	};
	
	\draw \foreach \a in {6,-6} \foreach \x in {0} {
		(\x:3.75)--(\x+\a:4.25) node {}
	};
	\draw[dotted]  \foreach \x in {0} {
		(\x-6:4.25)--(\x+6:4.25) 
	};
	
	\draw \foreach \x in {0} {
		%\draw[dotted] (\x:2.75)--(\x:3.5)  
		(\x:3.75)  node {}
		
	};

	\coordinate [label=left:$A_0$] (A) at (10.1:2.3); 
	\coordinate [label=right:$B_0$] (A) at (117.5:2.15); 
	\coordinate [label=right:$C_0$] (A) at (242.5:2.15); 
	\coordinate [label=left:$A_1$] (A) at (6.75:3); 
	%	\coordinate [label=right:$B_1$] (A) at (117.5:2.8); 
	%	\coordinate [label=right:$C_1$] (A) at (242.5:2.8); 
	\coordinate [label=left:$A_j$] (A) at (5.5:3.95); 
	%	\coordinate [label=right:$B_k$] (A) at (117.5:3.7); 
	%	\coordinate [label=right:$C_l$] (A) at (242.5:3.7); 

	\coordinate [label=center:$a$] (A) at (0:4.5); 
	\coordinate [label=center:$b$] (A) at (120:3); 
	\coordinate [label=center:$c$] (A) at (240:3); 
	
	\end{tikzpicture} \quad
	\centering
	\begin{tikzpicture}[thick,scale=1]%
	\draw \foreach \x in {0,90,180,270} {
		%(\x:3) node {} -- (\x+90,3))
		(\x:2) node {}  -- (\x+90:2)
		
	};
	\draw 	(0:2)--  (0:2.75)  node {};
	\draw[dotted] \foreach \x in {0} {
		(\x:2.75)--(\x:3.75)  
	};
	
	\draw[dotted] (175:2.5)--(185:2.5);
	\draw \foreach \a in {5,-5} \foreach \x in {0} {
		(\x:3.75)--(\x+\a:4.25) node {}
	};
	\draw \foreach \a in {10,-10} \foreach \x in {180} {
		(\x:2)--(\x+\a:2.5) node {}
	};
	\draw \foreach \x in {0} {
		%\draw[dotted] (\x:2.75)--(\x:3.5)  
		(\x:3.75)  node {}
		
	};
	
	\draw[dotted]  \foreach \x in {0} {
		(\x-5:4.25)--(\x+5:4.25) 
	};
	
	\coordinate [label=left:$A_0$] (A) at (2.5,0.35); 
	\coordinate [label=left:$B_0$] (A) at (0.85,2.1); 
	\coordinate [label=right:$C_0$] (A) at (-1.96,0); 
	\coordinate [label=left:$D_0$] (A) at (0.8,-2.1); 
	\coordinate [label=left:$A_1$] (A) at (3.15,0.35); 
	
	\coordinate [label=left:$A_j$] (A) at (3.85,0.3);

	\coordinate [label=center:$a$] (A) at (0:4.5); 
	\coordinate [label=left:$c$] (A) at (180:2.55);

	\end{tikzpicture}

	\caption{The graphs $G^3_{a,b,c,j}$ and $G^4_{a,c,j}$}
	\label{fig:G3abcj_G4acj}
\end{figure}
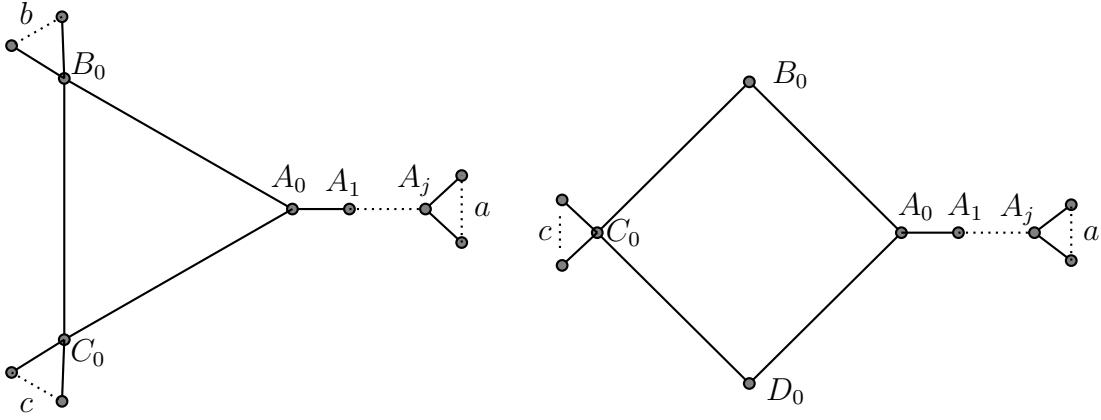

We conclude that the extremal graphs containing $C_3$ or $C_4$ are isomorphic to a graph of the form $G^3_{a,b,c,j}$ or $G^4_{a,c,j}$ with $a \in \{c,c+1\}$.
These graphs are shown in Figure~\ref{fig:G3abcj_G4acj}.

\section{Conditions on parameters of extremal  $G^3_{a,b,c,j}$ and $G^4_{a,c,j}$}\label{6}

%\subsection{$G$ contains $C_4$}
%
%We know already that $c \le a \le c+1$.

%Furthermore, when $a>0$ and $b=c=d=0$, we can find a better graph $G_3$ by connecting $B_0$ with $D_0$ and disconnecting $A_0$ with $B_0$.
%The matching number of the graphs are the same, while $W(G_3)-W(G)=j+a>0$.

%Furthermore, when $\lvert b-c \rvert \ge 2$ and $\min\{b,c\}>0$, the graph $G'$ with parameters $(b,c)$ replaced by $ \left( \lceil \frac{b+c}{2} \rceil , \lceil \frac{b+c}{2} \rceil \right)$ has a higher Wiener index and the same matching number,

%W(G)>W(G') since the essential part is $(b+c)^2+(b-1)(c-1)+1$

In this section, we prove the following proposition. 
\begin{prop}
	If the graph $G^3_{a,b,c,j}$ is an extremal graph, then $b=c=0$.
	If the graph $G^4_{a,c,j}$ (with $a\ge c$) is an extremal graph with $\max\{a,c,j\} \ge 1$, then $c\ge 1$.
\end{prop}

Take some graph $G=G^3_{a,b,c,j}$.

%If $b \ge c+2$ and $c>0$, then take $G'=G^3_{a,b',c',j}$ where $(b',c')= \left( \lceil \frac{b+c}{2} \rceil , \lfloor \frac{b+c}{2} \rfloor \right).$
%Now $m(G')=m(G)$ and $W(G')>W(G)$ by AM-GM, so $G$ is not extremal.

%Furthermore, when $\lvert b-c \rvert \ge 2$ and $\min\{b,c\}>0$, the graph $G'$ with parameters $(b,c)$ replaced by $ \left( \lceil \frac{b+c}{2} \rceil , \lceil \frac{b+c}{2} \rceil \right)$ has a higher Wiener index and the same matching number,

%W(G)>W(G') since the essential part is $(b+c)^2+(b-1)(c-1)+1$

If $a\ge b=1>c=0$, we choose the graph $G'=G^3_{a,0,0,j+1}.$
The graph $G$ was not an extremal graph, since $m(G')=m(G)$ and $W(G')-W(G)=j+a>0$.
If $a\ge b >1>c=0$, the graph $G'=G^4_{a,b-1,j}$ satisfies $m(G')=m(G)$ and $W(G')-W(G)=(b-1)(a+j+1)>0$.
When $a \ge b \ge c \ge 1$, the graph $G'=G^4_{a,b+c-1,j}$ satisfies $m(G')\le m(G)$ (equality when $2\nmid j$) and $W(G')-W(G)=(a+j+1)(c+b-1)-bc\ge 1$.
Hence $G$ was not an extremal graph, due to Proposition~\ref{monotonicity}.
So a graph $G^3_{a,b,c,j}$ with $(b,c)\not=(0,0)$ can not be extremal graph.
From now onwards we will write $G^3_{a,j}$ instead of $G^3_{a,0,0,j}$, this notation is shown in Figure~\ref{fig:G3aj_G4acj}.

Next, we prove that for an extremal graph $G^4_{a,c,j}$ with more than $4$ vertices, we have that $c>0$.
Take an extremal graph $G^4_{a,c,j}$
Since an extremal graph satisfies $\lvert a -c \rvert \le 1$ as shown in Subsection~\ref{sec:C4}, we can take $a=c=0$. (since $G^4_{1,0,j}=G^4_{0,0,j+1}$)

Observe that $m(G^4_{0,0,j})=m(G^3_{0,j+1})$ and $W(G^4_{0,0,j})<W(G^3_{0,j+1})$, from which the conclusion follows.

    	\section{Calculating $W(\U(n,m))$}\label{7}

From the results of previous sections, the extremal graphs are of the form $G^3_{a,j}$, $G^4_{a,a-1,j}$ with $a\ge 2$ or $G^4_{a,a,j}$ with $a \ge 1$.
We will determine the extremal values $W(\U(n,m))$ for $m\ge 2$.

When $G^3_{a,j}$ is a graph with $n$ vertices and matching number $m$,
then $n=a+j+3$ and $m=2+\lfloor \frac j2 \rfloor$.
So $j \in \{2m-4,2m-3\}$ and $a=n-(j+3).$
Notice that $W\left (G^3_{n-2m+1,2m-4}\right )\le W\left (G^3_{n-2m,2m-3}\right )$ with equality iff $n=2m$, in which case $a=0$ and so we actually look only to the same graph.

%$W\left (G^3_{n-2m,2m-3}\right )=2-\frac83 m^3+2m^2+\frac53 m+2m^2n-3mn-2n+n^2$

Next, note that $m(G^4_{a,a-1,j})=m(G^4_{a,a,j})=2+\lfloor \frac {j+1}2 \rfloor$ and $n(G^4_{a,a-1,j})+1=n(G^4_{a,a,j})=j+4+2a.$
So $j \in \{2m-5,2m-4\}$. 

According to the parity of $n$, we can compare the Wiener index of the corresponding graphs $G^4_{a,a-1,j}$ and $G^4_{a,a,j}$.
We check that 
$$W\left(G^4_{n/2-m,n/2-m,2m-4}\right)-W\left(G^4_{n/2-m+1,n/2-m,2m-5}\right)=\frac{1}{4}(n-2m)(n+2m-4)\ge 0$$
where equality cannot occur, since  $m \ge 2$ and $n \ge 2m+2$ as we need $n/2-m\ge 1$.

When $n$ is odd, we have that
\begin{align*}
W\left(G^4_{n/2-m+1/2,n/2-m-1/2,2m-4}\right)-W\left(G^4_{n/2-m+1/2,n/2-m+1/2,2m-5}\right)\\
=\frac14 (n-2m)(n+2m-4)+m-\frac{13}4\ge 0
\end{align*}
since $n\ge 2m+1$ and $m\ge 2$. Equality occurs only when $n=5$ and $m=2$, but then $2m-5<0$, implying that the graph $G^4_{n/2-m+1/2,n/2-m+1/2,2m-5}$ does not exist. 

%Remark that it is easy to see that we will take $j$ as large as possible.

To finish the search for the extremal graph, we have to compare the graphs of the form $G^4$ with the $G^3$ graph.
When $n$ is odd, 
\begin{align*}
&W\left(G^4_{n/2-m+1/2,n/2-m-1/2,2m-4}\right)-W\left (G^3_{n-2m,2m-3}\right )\\&=n+2m^3+nm+\frac12 mn^2-\frac12n^2-\frac72m-2nm^2+\frac52.
\end{align*}
Taking $n=2m+k$, we get that this value equals $1+\frac12 (k-3)(k+1)(m-1)$ which is strictly positive for $k\ge 3$ and strictly negative for $k=1$.
When $n$ is even,
$$W\left(G^4_{n/2-m,n/2-m,2m-4}\right)-W\left (G^3_{n-2m,2m-3}\right )=4+n+2m^3+nm+\frac12 mn^2-\frac12n^2-4m-2nm^2.$$
Taking $n=2m+k$, this expression equals $2+\frac12(k^2-2k-4)(m-1)$.
This is smaller or equal than zero for $k\in \{0,2\}$ with equality iff $m=2$.
When $k\ge 4$, it is strictly positive.
We summarize all those results in the following theorem with Figure~\ref{fig:G3aj_G4acj} showing the exact picture of the extremal graphs.

%\subsection{Case $G^4_{a,a-1,j}$ and $n$ odd}
%
%We have that $j=2m-4$ and $a=\frac{n-2m+1}{2}$.
%
%
%\subsection{Case $G^4_{a,a-1,j}$ and $n$ even}
%We have that $j=2m-5$ and $a=\frac{n-2m}{2}+1$.
%
%\subsection{Case $G^4_{a,a,j}$ and $n$ odd}
%We have that $j=2m-5$ and $a=\frac{n-2m+1}{2}$.
%
%\subsection{Case $G^4_{a,a,j}$ and $n$ even}
%We have that $j=2m-4$ and $a=\frac{n-2m}{2}$.

%We conclude that 

%\[ W(\U(n,m))=\begin{cases} 
%2-\frac83 m^3+2m^2+\frac53 m+2m^2n-3mn-2n+n^2 & \mbox{ if } n \le 2m+2 \\
%6-n-\frac23 m^3+\frac12 n^2-2nm+2m^2+\frac12 mn^2-\frac73m & \mbox{ if } n \ge 2m+4 \mbox{ and } n \mbox{ is even} \\
%\frac92-n-\frac23 m^3+\frac12 n^2-2nm+2m^2+\frac12 mn^2-\frac{11}6m & \mbox{ if } n \ge 2m+3 \mbox{ and } n \mbox{ is odd}
%\end{cases}
%\]

\begin{thr}\label{mainthr}
	Let $G \in \U(n,m)$, where $2 \le m \le \lfloor \frac n2 \rfloor$.
	
	\begin{itemize}
		\item If  $n \le 2m+2$, then $W(G) \le 2-\frac83 m^3+2m^2+\frac53 m+2m^2n-3mn-2n+n^2$ with equality iff  $G=G^4_{0,0,0},G^3_{0,1}$ for $(n,m)=(4,2)$, $G=G^4_{1,1,0},G^3_{2,1}$ for $(n,m)=(6,2)$ and $G=G^3_{n-2m, 2m-3}$ otherwise.
		\item If $n \ge 2m+3$ and $n$ is odd, then $W(G) \le \frac92-n-\frac23 m^3+\frac12 n^2-2nm+2m^2+\frac12mn^2-\frac{11}6m$ with equality iff $G=G^4_{n/2-m+1/2,n/2-m-1/2,2m-4}$.
		\item If $n \ge 2m+4$ and $n$ is even, then $W(G) \le 6-n-\frac23 m^3+\frac12 n^2-2nm+2m^2+\frac12 mn^2-\frac73m$ with equality iff $G=G^4_{n/2-m,n/2-m,2m-4}$.
	\end{itemize}	
\end{thr}

    	\subsection*{Acknowledgment}
The author is very grateful to Jan Bok for presenting this problem in Ghent at GGTW $2017.$
%to the anonymous referees for the time they dedicated to this article and for the comments they made which improved the first version of this manuscript.

%%%%%%%%%%%% BIBLIO %%%%%%%%%%%%%%%%%%%%


\begin{thebibliography}{99}
	
	\baselineskip=0.0in
	
	\bibitem{D}
	P. Dankelmann, Average distance and independence number, {\it Discr. Appl. Math.\/} {\bf	51} (1994) 75-83.
	

	
	\bibitem{DZ}
	Z. Du, B. Zhou, Minimum Wiener indices of trees and unicyclic graphs of given matching number, {\it MATCH Commun. Math. Comput. Chem.\/} {\bf 63} (2010) 101--112.
	
	\bibitem{KST}
	M. Knor, R. \v{S}krekovski, A. Tepeh, Mathematical aspects 
	of Wiener index, {\it Ars Math. Contemp.\/} {\bf 11}
	(2016) 327--352.
	
 
	
	
	

	
	
\end{thebibliography}
\end{document}